\documentclass{article}
\usepackage{amssymb,amsmath}
\usepackage{anysize}
\usepackage{psfrag}
\newtheorem{thm}{Theorem}[section]
\newtheorem{def.}{Definition}[section]
\newtheorem{rem}{Remark}[section]
\newtheorem{prop}{Proposition}[section]

\newtheorem{lem}{Lemma}[section]

\numberwithin{table}{section}

\begin{document}

\title{Determinants of Rational Knots}
        \author{Louis H. Kauffman\\
        Department of Mathematics, Statistics and Computer Science\\
        University of Illinois at Chicago\\
        851 S. Morgan St., Chicago IL 60607-7045\\
        USA\\
        \texttt{kauffman@uic.edu}\\
        and\\
        Pedro Lopes\\
        Department of Mathematics\\
        Instituto Superior T\'ecnico\\
        Technical University of Lisbon\\
        Av. Rovisco Pais\\
        1049-001 Lisbon\\
        Portugal\\
        \texttt{pelopes@math.ist.utl.pt}\\
}
\date{July 17, 2009}
\maketitle

\begin{abstract}
We study the Fox coloring invariants of rational knots. We express the propagation of the colors down the twists of these knots and ultimately the determinant of them with the help of finite increasing sequences whose terms of even order are even and whose terms of odd order are odd.
\end{abstract}

\bigbreak


\bigbreak

2000MSC: 05A99, 57M27

\section{Introduction}

\noindent

\begin{figure}[h!]
    \psfrag{a}{\huge $a$}
    \psfrag{b}{\huge $b$}
    \psfrag{c}{\huge $c$}
    \psfrag{2b=a+c}{\huge $2b=a+c$}
    \centerline{\scalebox{.50}{\includegraphics{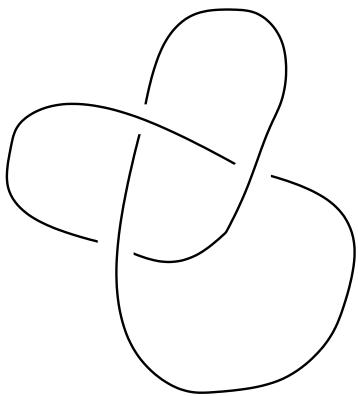}}}
    \caption{A diagram of the trefoil knot.}\label{Fi:trefdiagram}
\end{figure}

A knot is an embedding of the circle into three-dimensional space. Knots that are obtained from one another by continuous deformation of their embeddings are said to be equivalent. The classification of these equivalence classes is still an open problem. In order to study a knot one usually resorts to projecting it into a plane in such a way that the singularities of the projection look locally like the crossing of two line segments. At these crossings the line that goes under in the embedding is broken in the projection giving rise to the so-called knot {\it diagram} (see Figure \ref{Fi:trefdiagram}).
\bigbreak

Knot diagrams are planar networks of arcs. The Reidemeister moves (\cite{lhKauffman}) are defined to be local transformations on the diagrams. The Reidemeister theorem states that, given two diagrams, the knots corresponding to these diagrams are  equivalent if and only if  these diagrams are obtained from one another by a finite number of the Reidemeister moves. Whenever a mathematical object is associated to a knot diagram in such a way that this object remains the same when Reidemeister moves are performed on the diagram, then this object is called a knot invariant. By Reidemeister's theorem, a knot invariant is an invariant of the equivalent class of the knot. The knot invariant we will be concerned with in this article is called {\it Fox coloring} (\cite{lhKauffman}).

\bigbreak
We now elaborate on Fox colorings. Given any
knot diagram of the knot under study, the arcs of the diagram stand for algebraic unknowns and the equation ``the sum of the under-arcs equals twice the over-arc'' is read off at each crossing (see Figure \ref{Fi:crosseq}). The
system of equations so obtained has a matrix of coefficients over
the integers. The equivalence class of this matrix under
elementary operations on  matrices (as listed on page 50 of \cite{Lickorish}) over the integers constitutes an invariant of the
knot under study. The elementary operations on matrices over the integers are standard row and column operations plus enlarging the matrix by a row and column of zeros except for the diagonal entry of $1$, or removing such a row and column.

Other invariants of the knot stemming from this one are the following.
The invariant factors of the equivalence class of the matrix; any first minor (also known as the
determinant of the knot); the number of solutions of the system of
equations over the integers modulo $r$ (also known as the number of Fox $r$-colorings). The minimum
number of distinct colors it takes to represent such a nontrivial
solution (i.e., a solution such that at least two unknowns
take on different values) across the diagrams that represent this knot, is a subtle and interesting invariant (\cite{klopes}).

Some of these invariants are easy to calculate and, at the same
time, seem to be effective enough to be interesting (number of
colorings, see \cite{pLopes}, \cite{DL}). Others are hard to determine and
related to elusive conjectures (minima of colors, see \cite{hk},
\cite{kl-t}, \cite{aps}, \cite{klopes}).

In \cite{DL} a brute force approach was used: a knot diagram was
introduced into a computer program that converted it into a
coloring system of equations. Then, candidates for solutions of the
coloring system of equations were tested. Finally, the actual
solutions were counted (\cite{DL}).

In \cite{klopes} the focus was on a
particular class of knots and a general formula was obtained for
the number of colorings. For this class of knots, a
Gaussian reduction of the systems of equations was developed by
inspection of the knot diagram under study. Some unknowns
were enough to express the rest of them. This led to the
development of a black-box approach to diagrams. What
is inside a given part of the diagram, the so-called
black-box, is ignored. All that matters is that there is a
``color input'' and a ``color output'' to this part of the diagram. By calculating the propagation of
the ``color input'' down the black-box, the ``color output'' is specified in terms of the ``color
input'' and the ``size'' of the black-box. Finally we equate this ``color output''  to the ``color input'' in order to obtain the coloring system of  equations and other invariants like the determinant of the knot (we elucidate this below, see Section \ref{sect : black-box}).
This technique can then be applied to more complicated situations.
Whenever a knot diagram involves some combination of these black-boxes,
 we use this approach, thus reducing the number of
unknowns and the number of equations we have to work with.

A particular situation comes up with a ``twist'', which is obtained by twisting two line segments an assigned number of times (see the left-hand side of Figure \ref{Fi:tref} for a particular twist).
Since the propagation of colors down a twist is
well-understood (see \cite{klopes}) and since rational knots (see Section \ref{sect : rational} below) can be
regarded as special stackings of twists , we look into
 calculating the propagation of colors down the twists of a rational knot. We denote $R(n\sb{1}, n\sb{2}, \dots , n\sb{N})$ the rational knot with $n\sb{i}$ crossings on the $i$-th twist. The formulas we obtain involve three polynomials in the number
of crossings on each twist, the $n\sb{i}$'s. These polynomials are closely related to the set of increasing sequences of integers from $\{1, 2, \dots , N\}$ whose terms of even order are even and whose terms of odd order are odd, call this set $IEO[N]$. The coefficients of these polynomials are all equal to $1$. In each monomial the power of each $n\sb{i}$ is either $1$ or $0$. One of the polynomials, call it $p\sb{N}$, has as many monomials as there are sequences in $IEO[N]$. As for the other two, one of them, call it $p\sb{N}\sp{e}$, is obtained from  $p\sb{N}$ by deleting the monomials which are products of an odd number of $n\sb{i}$'s. The other polynomial, call it $p\sb{N}\sp{o}$, equals $p\sb{N}-p\sb{N}\sp{e}$. Moreover, the relation between the monomials of $p\sb{N}$ and the sequences in $IEO[N]$ is as follows. If $(u\sb{1}, u\sb{2}, \dots , u\sb{k})\in IEO[N], \,  (1\leq k\leq N)$, then $n\sb{u\sb{1}} n\sb{u\sb{2}} \dots n\sb{u\sb{k}}$ is a monomial in $p\sb{N}$ (and conversely).

We prove the following result in Section \ref{sect:colratlinkgen}.

\bigbreak

\begin{thm}\label{thm:thethm} Given a non-negative integer $I$, consider $2I+1$ integers $n\sb{1}, n\sb{2}, \dots , n\sb{2I},
n\sb{2I+1}$.
\begin{enumerate}
\item The coloring system of equations of $R(n\sb{1}, n\sb{2},
\dots , n\sb{2I})$ and its determinant are, respectively
\[
(b-a)p\sb{2I}\sp{e} = 0 \qquad  \text { and } \qquad \det
R(n\sb{1}, n\sb{2}, \dots , n\sb{2I}) = p\sb{2I}\sp{e}.
\]

\item The coloring system of equations of $R(n\sb{1}, n\sb{2},
\dots , n\sb{2I}, n\sb{2I+1})$ and its determinant are,
respectively
\[
(b-a)p\sb{2I+1}\sp{o} = 0 \qquad  \text { and } \qquad \det
R(n\sb{1}, n\sb{2}, \dots , n\sb{2I}, n\sb{2I+1}) =
p\sb{2I+1}\sp{o}.
\]
\end{enumerate}
\end{thm}

\bigbreak

This article is organized as follows. In Section \ref{sect :
black-box} we introduce the background material and the
notion of black-box in the context of knot diagrams and colorings. In Section \ref{sect : rational} we
introduce the diagrammatics of rational knots. In Section \ref{sect:colratlink} we
calculate the coloring equation and the determinant of a few
particular rational knots and single out patterns in their expressions which involve the polynomials referred to above. The polynomials  are formally introduced
in Section \ref{sect : polynomials} and relations among them are proved. In Section~\ref{sect:colratlinkgen} we prove  Theorem \ref{thm:thethm} stating the form of the coloring equation and the determinant of a rational knot $R(n\sb{1}, \dots , n\sb{N})$. In Section \ref{sect:count} we give an alternative approach to the determinants of rational  knots and links by counting the spanning trees in their checkerboard graphs. This leads to a simple recursion formula for these determinants. The outputs of these formulas are identical to the results of  Theorem \ref{thm:thethm}

\bigbreak

\subsection{Acknowledgements} \label{subsect:ack}

\noindent

P.L. acknowledges support by {\em Programa Operacional
``Ci\^{e}ncia, Tecnologia, Inova\c{c}\~{a}o''} (POCTI) of the {\em
Funda\c{c}\~{a}o para a Ci\^{e}ncia e a Tecnologia} (FCT)
cofinanced by the European Community fund FEDER. He also thanks
the staff at IMPA and especially his host, Marcelo Viana, for
hospitality during his stay at this Institution.

We thank the referee for useful suggestions which improved this article.

\bigbreak

\section{The Black-box Approach}\label{sect : black-box}

\noindent

We introduce some of the objects we will be dealing with in this
article.

\begin{figure}[h!]
    \psfrag{3}{\Huge $3$}
    \psfrag{a}{\huge $a$}
    \psfrag{b}{\huge $b$}
    \psfrag{2b-a}{\huge $2b-a$}
    \psfrag{3b-2a}{\huge $3b-2a$}
    \psfrag{b+3(b-a)}{\huge $b+3(b-a)$}
    \psfrag{a+3(b-a)}{\huge $a+3(b-a)$}
    \centerline{\scalebox{.50}{\includegraphics{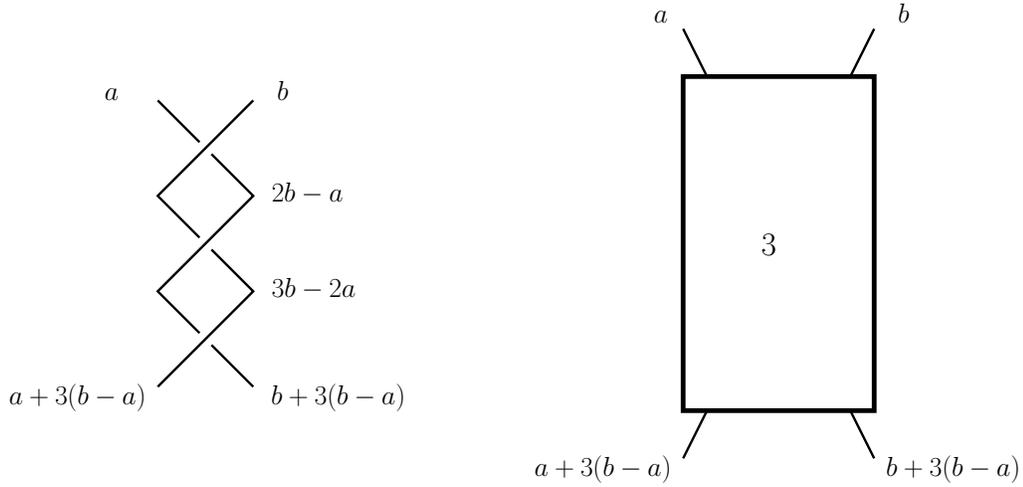}}}
    \caption{The braid $\sigma\sb{1}\sp{3}$ whose closure
yields a diagram of the trefoil, and setting up a Fox coloring.
 Left: the integral view. Right: the black box
approach.}\label{Fi:tref}
\end{figure}

\begin{def.}[Twist]\label{def : twist} A
 twist is a portion of a diagram which is
given by $\sigma\sb{1}\sp{n}$ where $\sigma\sb{1}$ is a standard
generator of the braid group on two strands, $B\sb{2}$ $(\text{see
} \cite{Birman})$. The left hand-side of Figure \ref{Fi:tref}
depicts $\sigma\sb{1}\sp{3}$.
\end{def.}

\begin{def.}[Fox Coloring, \cite{Fox}]\label{def : foxc} A Fox coloring of a
knot or tangle diagram is the assignment of integers to the arcs
of the diagram such that at each crossing twice the integer
assigned to the over-arc equals the sum of the integers assigned
to the under-arcs meeting at this crossing, see Figure \ref{Fi:crosseq}.
\begin{figure}[h!]
    \psfrag{a}{\huge $a$}
    \psfrag{b}{\huge $b$}
    \psfrag{c}{\huge $c$}
    \psfrag{2b=a+c}{\huge $2b=a+c$}
    \centerline{\scalebox{.50}{\includegraphics{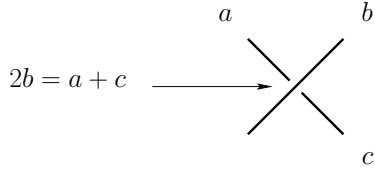}}}
    \caption{The equation associated to the unknowns at a crossing.}\label{Fi:crosseq}
\end{figure}We equivalently express this
equality as ``the sum of the under-arcs equals twice the over-arc''.

Furthermore, these equations and their solutions will be considered modulo $r$ for
given integers $r>1$.

In the present work, we do not specify this modulus
so that the formulas we obtain are as general as possible.
Moreover, in the literature, the modulus is often taken to be a
prime number. Here we allow the modulus to take on any
composite integer value greater than $2$, $r\in \mathbb{Z}\sp{+}
\setminus \{ 1, 2 \}$.
\end{def.}

\begin{def.}[Coloring System of Equations]\label{def : colsys}
A Coloring System of Equations $($CSE$)$ of a knot diagram is the system of
equations assigned to the knot diagram by regarding the arcs as
unknowns and reading the equation $2b=a+c$ at each crossing $($see Figure \ref{Fi:crosseq}$)$. Usually we solve these CSE's over the integers $\mod r$, for some modulus $r$. The solutions of this CSE are the Fox colorings of the knot diagram with
respect to the specified modulus. The Coloring Matrix $($CM$)$ of the knot
diagram is the matrix of the coefficients of the CSE of this knot diagram. This is a square matrix since
any diagram of a non-trivial knot has as many arcs as crossings.
Any matrix obtained from the CM by elementary matrix
operations - \cite{Lickorish}, page 50 - is also known as a coloring matrix .
The system of equations associated to this new matrix is also
a coloring system of equations for this knot.
\end{def.}


\begin{prop} Consider a modulus and a knot $K$. The number of Fox
colorings in this modulus for any knot diagram of $K$ is an
invariant of the knot. Also, the absolute value of any first minor of the
coloring matrix is an invariant of
the knot known as the ``determinant of the knot'' and denoted $\det K$.

Non-trivial colorings occur only for those moduli which are not relatively prime to the determinant of the knot.
\end{prop} Proof: See \cite{lhKauffman}. $\hfill \blacksquare $

\bigbreak

Consider the trefoil knot regarded as the closure of
$\sigma\sb{1}\sp{3}(\in B\sb{2})$ (see left hand-side of Figure
\ref{Fi:tref}). Endowing this diagram with a Fox coloring amounts
to assigning integers to the arcs of the diagram and making sure
that, at each crossing, the rule ``$2b=a+c$'' holds (in a given modulus). On the other
hand, working our way down from the top of the braid in the left
hand-side of Figure \ref{Fi:tref}, we can regard each crossing as
a rule that sets the emerging under-arc equal to twice the
over-arc minus the in-coming under-arc.

\begin{figure}[h!]
    \psfrag{...}{\Huge $\vdots $}
    \psfrag{n}{\Huge $n$}
    \psfrag{nth}{\Large $n$-th xing}
    \psfrag{ith}{\Large $i$-th xing}
    \psfrag{(i+1)b-ia}{\huge $(i+1)b-ia$}
    \psfrag{a}{\huge $a$}
    \psfrag{b}{\huge $b$}
    \psfrag{2b-a}{\huge $2b-a$}
    \psfrag{3b-2a}{\huge $3b-2a$}
    \psfrag{b+n(b-a)}{\huge $b+n(b-a)$}
    \psfrag{a+n(b-a)}{\huge $a+n(b-a)$}
    \centerline{\scalebox{.50}{\includegraphics{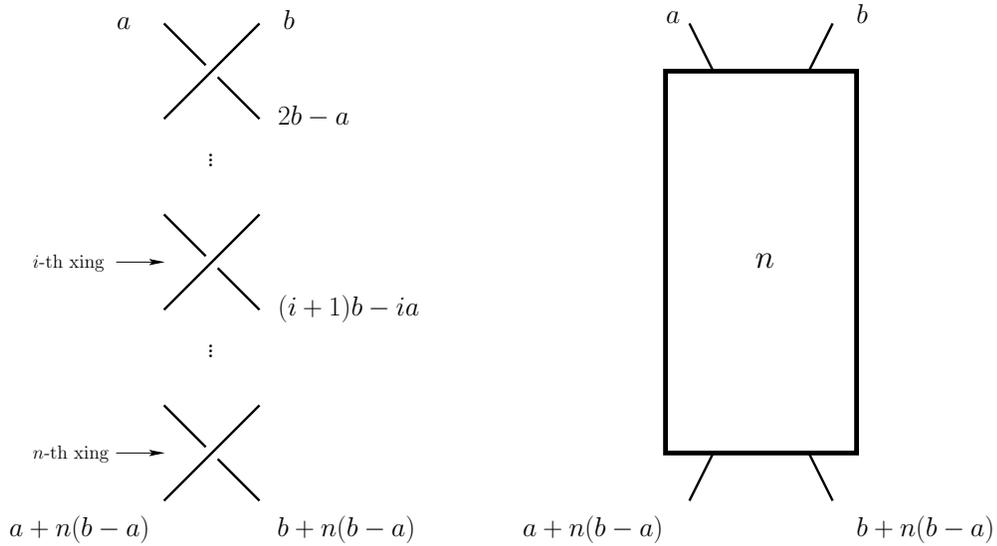}}}
    \caption{The braid $\sigma\sb{1}\sp{n}$ whose closure
yields a diagram of the $T(2, n)$ torus knot, and endowed with a
Fox coloring. Left: the integral view. Right: the ``black box''
approach.}\label{Fi:sigman}
\end{figure}

We believe the assignments to arcs in Figures \ref{Fi:tref} and
\ref{Fi:sigman} to be self-explanatory, Figure \ref{Fi:tref} for
the $n=3$ case and Figure \ref{Fi:sigman} for the general $n$
case (see also \cite{klopes}). Eventually the bottom arcs are reached and in order to obtain
a coloring from this assignment, the color obtained at the bottom
right arc has to equal $b$ and the color at the bottom left arc
has to equal $a$. These two equations are equivalent to each other and to:

\begin{equation*}
\begin{cases}
n(b-a) = 0 \\
0 = 0.
\end{cases}
\end{equation*}

The determinant of the trefoil is then $3$, which is the $n=3$ case. We kept here the $0=0$
equation for consistency. In the sequel we will drop it whenever
we arrive at a situation where two equations in a system are
equivalent.

It was not the specifics at each crossing of the braid that
allowed us to write down the coloring system of equations (and
consequently the determinant of the knot) but the fact that the
colors at the bottom strands are expressed in a linear fashion in
terms of the colors at the top strands and the number of
crossings. This led us to the black-box approach which is
depicted on the right-hand side of Figure \ref{Fi:tref} for the $3$-crossing instance. Here we
regard the colors $a$ and $b$ at the top strands as the color
input to a black-box labelled with a $3$ which transforms them
into a color output at the bottom strands, $a+3(b-a)$ and
$b+3(b-a)$, from left to right. Figure \ref{Fi:sigman} addresses the case with $n$ crossings.

\bigbreak

The following Proposition summarizes the situation so far and includes a
formula for the calculation of the number of colorings for a $T(2, n)$
knot. A proof can be found in \cite{klopes}.

\begin{prop}\label{prop:colortwist} Given a non-zero integer $n$ and integers $a$, $b$,
and $r(>1)$, consider $\sigma\sb{1}\sp{n}(\in B\sb{2})$.  Assume
$a$ and $b$ is the color input on the top of the braid $($see
Figure \ref{Fi:sigman}$)$. Then the arc emerging from the
$n$-th crossing receives color $b+n(b-a)$. The
coloring system of equations of the $T(2, n)$ torus knot is formed by $n(b-a) =
0$ mod $r$ $($and the trivial equation$)$. It has $\gcd(n, r) r$
solutions. Moreover, $\det T(2, n) = n$.
\end{prop}

\bigbreak

Furthermore, and in view of the preceding discussion, we set forth the following Definition.

\begin{def.}
Consider the right-hand side of Figure \ref{Fi:sigman}. It is
formed by a rectangular box with two strands on top of it and two
strands on the bottom. The box has a label $n$ inside. The top
strands are labeled $a$ and $b$ from left to right and the bottom
strands are labeled $a+n(b-a)$ and $b+n(b-a)$. This ensemble
stands for a twist equipped with a Fox-coloring induced by
the input $a$ and $b$ on the top. We denote this the Black Box
$($for the twist$)$ when we ignore the color input. Otherwise
we denote it the Black Box $($for the twist$)$ endowed
with a Fox coloring.
\end{def.}

\bigbreak

\section{Rational Knots and Their Representation}\label{sect : rational}

\noindent

In this section we describe the diagrammatic representation of rational knots that fits our needs herein. A
formal introduction to this material may be found in \cite{kl-k}
and \cite{kl-t}.

\subsection{Rational Knots as Special Closures of $4$-Strand Braids.}

\noindent

A rational knot can be regarded as a special composition of twists using the generators $\sigma\sb{1}$ and $\sigma\sb{2}$ as well as their inverses, from the braid group on four strands $B\sb{4}$.

\begin{figure}[h!]
    \psfrag{n1}{\Huge $n\sb{1}$}
    \psfrag{n2}{\Huge $n\sb{2}$}
    \psfrag{n3}{\Huge $n\sb{3}$}
    \psfrag{n4}{\Huge $n\sb{4}$}
    \psfrag{n5}{\Huge $n\sb{5}$}
    \centerline{\scalebox{.20}{\includegraphics{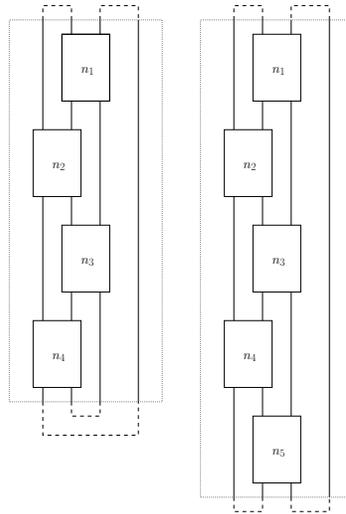}}}
    \caption{Rational knots on $4$ and $5$ twists.}\label{Fi:newclbraidrational}
\end{figure}

For each positive integer $N$, we choose $N$ integers $n\sb{1}, n\sb{2}, \dots , n\sb{N}$ and form the braid:
\[
\sigma\sb{2}\sp{n\sb{1}}\sigma\sb{1}\sp{n\sb{2}} \dots  \sigma\sb{2}\sp{n\sb{N-1}}\sigma\sb{1}\sp{n\sb{N}}, \qquad \qquad \text{ for even } N,
\]
and
\[
\sigma\sb{2}\sp{n\sb{1}}\sigma\sb{1}\sp{n\sb{2}} \dots  \sigma\sb{1}\sp{n\sb{N-1}}\sigma\sb{2}\sp{n\sb{N}}, \qquad \qquad \text{ for odd } N.
\]
See Figure \ref{Fi:newclbraidrational} for illustrating examples with $N=4$ and $N=5$ (ignore the broken lines for the moment). We elucidate the sign convention below.

Finally, a {\bf special closure} of these braids is performed in order to obtain the corresponding rational knot which we have been denoting $R(n\sb{1}, n\sb{2}, \dots , n\sb{N})$. These plat closures are depicted by the broken lines in Figure \ref{Fi:newclbraidrational}. We believe Figure \ref{Fi:newclbraidrational} is now self-explanatory for the $N=4$ and $N=5$ cases, elucidating the general case. We note that the special closure in the bottom for the even $N$ case (called {\bf denominator closure}) is different from the odd $N$ case (called {\bf numerator closure}).

\bigbreak

\subsection{Checkerboard shadings and sign convention for
crossings}\label{subsect : checkershades}

\noindent

The signs at crossings we define for diagrams of rational tangles
are consistent with the checkerboard shading of knot diagrams in
the way defined below.

\begin{figure}[h!]
    \psfrag{+}{\huge $(+)$}
    \psfrag{m}{\huge $(-)$}
    \centerline{\scalebox{.50}{\includegraphics{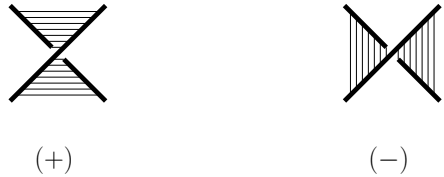}}}
    \caption{The signs of the crossings induced by the checkerboard
shading.}\label{Fi:cross}
\end{figure}

\begin{figure}[h!]
    \psfrag{r4-3}{\huge $\mathbf{R(4, -3)}$}
    \psfrag{r43}{\huge $\mathbf{R(4, 3)}$}
    \centerline{\scalebox{.40}{\includegraphics{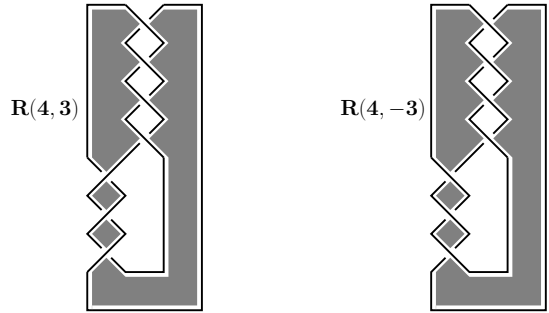}}}
    \caption{Sign convention and checkerboard shadings.}\label{Fi:rationalexgf}
\end{figure}

\begin{def.}[Checkerboard shadings and sign conventions]\label{def:checkandsign}
A checkerboard shading of a knot diagram is the shading of some
regions of the knot diagram such that the four regions meeting at
each crossing have the following property. Whenever two of these
four regions share a common boundary (which is an arc of the diagram), then one of them is shaded
and the other one is not. The signs at crossings are then induced
by the checkerboard shading in the following way. Consider a
crossing of the diagram along with its over-arc. Rotate this
over-arc counterclockwise about the crossing. If the region swept
out by the over-arc is the shaded region then the crossing is
positive. Otherwise it is negative $($see Figure
\ref{Fi:cross} and \ref{Fi:rationalexgf}$)$.
\end{def.}

\bigbreak

 As
a rule, we do not shade the unbounded region around the diagram.
See Figure \ref{Fi:rationalexgf} for the illustrative examples of
the sign convention and its connection to the checkerboard
shadings.

\bigbreak

\section{The Coloring System of Equations of a Rational Knot: Particular Examples}\label{sect:colratlink}

\noindent

In this section we will be always referring to a rational knot
with $N$ given twists ($1\leq N \leq 5$) where the $i$-th twist has $n\sb{i}$
crossings, $R(n\sb{1}, n\sb{2}, \dots , n\sb{N})$. The signs of the $n\sb{i}$'s are induced by the
checkerboard shading of the diagram where the unbounded exterior
of the diagram is not shaded, in the way described by Definition \ref{def:checkandsign}. The general look of the diagrams of
these knots that we will be using in the sequel, materializes into
the left- and right-hand sides of Figure \ref{Fi:newclbraidrational}
for $N=4$ and $N=5$, respectively. We
usually do not draw the rightmost strand for the rational knots nor do we perform all the plat
closures (as we did in Figure \ref{Fi:newclbraidrational}), in order to graphically bring out the
colors the arcs at the bottom are receiving. Nonetheless, the
relevant closure of $R(n\sb{1}, n\sb{2}, \dots , n\sb{N})$ is to
be assumed. In particular, in Figure \ref{Fi:rn1n2} the arc assigned $m\sb{2}$ should be connected to the arc assigned $r\sb{2}$
and the arc assigned $l\sb{2}$ should be connected to the arc assigned $a$, following the pattern established in Figure \ref{Fi:newclbraidrational} (even case).

\subsection{The rational knots $R(n\sb{1})$
through $R(n\sb{1}, n\sb{2}, n\sb{3}, n\sb{4}, n\sb{5}).$}\label{ss:particularcases}

\noindent

In the next subsections we work out the coloring system of
equations and the determinant of the indicated rational knots.

\subsubsection{$R(n\sb{1})$.}\label{sss:rn1}

\noindent

$R(n\sb{1})$ is the torus knot $T(2, n\sb{1})$. As was seen
above (Proposition \ref{prop:colortwist}), the coloring
system of equations reduces to:
\[
(b-a)n\sb{1} = 0
\]
and the determinant is

\[
\det R(n\sb{1}) = n\sb{1}.
\]

\subsubsection{$R(n\sb{1}, n\sb{2})$} \label{sss:rn1n2}

\noindent

\begin{figure}[h!]
    \psfrag{a}{\Huge $a$}
    \psfrag{b}{\Huge $b$}
    \psfrag{l1}{\Huge $l\sb{1}$}
    \psfrag{m1}{\Huge $m\sb{1}$}
    \psfrag{r1}{\Huge $r\sb{1}$}
    \psfrag{l2}{\Huge $l\sb{2}$}
    \psfrag{m2}{\Huge $m\sb{2}$}
    \psfrag{r2}{\Huge $r\sb{2}$}
    \psfrag{n1}{\Huge $n\sb{1}$}
    \psfrag{n2}{\Huge $n\sb{2}$}
    \centerline{\scalebox{.25}{\includegraphics{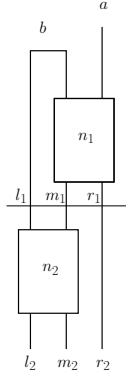}}}
    \caption{$R(n\sb{1}, n\sb{2})$ upon denominator closure.}\label{Fi:rn1n2}
\end{figure}

See Figure \ref{Fi:rn1n2}.

Using the black-box approach developed above, we have, at the
cross-section $(l\sb{1}, m\sb{1}, r\sb{1})$:
\begin{align}\notag
l\sb{1}&=b , \\ \notag m\sb{1}&=b+(b-a)n\sb{1}, \\ \notag
r\sb{1}&=a+(b-a)n\sb{1},
\end{align}
and so at the $(l\sb{2}, m\sb{2}, r\sb{2})$ cross-section we have:
\begin{align}\notag
l\sb{2}&=l\sb{1}+n\sb{2}(m\sb{1}-l\sb{1}) = b+(b-a)n\sb{1}n\sb{2} = a+(b-a)( 1 + n\sb{1}n\sb{2}) , \\
\notag m\sb{2}&=m\sb{1}+n\sb{2}(m\sb{1}-l\sb{1}) = b+(b-a)(n\sb{1}+n\sb{1}n\sb{2}) = a+(b-a)( 1 + n\sb{1}+n\sb{1}n\sb{2}) , \\
\notag r\sb{2}&=r\sb{1}=a+(b-a)n\sb{1}.
\end{align}

At this point, we perform a special closure of this
rational tangle for even $N$ i.e., we identify the arcs that receive colors
$m\sb{2}$ and $l\sb{2}$, and the arcs that receive colors
$a$ and $l\sb{2}$. From setting the corresponding pairs of colors equal we
obtain two equations which are equivalent to:
\[
(b-a)(1+n\sb{1}n\sb{2})= 0
\]
and thus the determinant of $R(n\sb{1}, n\sb{2})$ is
\[
\det R(n\sb{1}, n\sb{2}) = 1+n\sb{1}n\sb{2}.
\]

\subsubsection{$R(n\sb{1}, n\sb{2}, n\sb{3})$}\label{sss:rn1n2n3}

\noindent

\begin{figure}[h!]
    \psfrag{a}{\Huge $a$}
    \psfrag{b}{\Huge $b$}
    \psfrag{l1}{\Huge $l\sb{1}$}
    \psfrag{m1}{\Huge $m\sb{1}$}
    \psfrag{r1}{\Huge $r\sb{1}$}
    \psfrag{l2}{\Huge $l\sb{2}$}
    \psfrag{m2}{\Huge $m\sb{2}$}
    \psfrag{r2}{\Huge $r\sb{2}$}
    \psfrag{l3}{\Huge $l\sb{3}$}
    \psfrag{m3}{\Huge $m\sb{3}$}
    \psfrag{r3}{\Huge $r\sb{3}$}
    \psfrag{n1}{\Huge $n\sb{1}$}
    \psfrag{n2}{\Huge $n\sb{2}$}
    \psfrag{n3}{\Huge $n\sb{3}$}
    \centerline{\scalebox{.250}{\includegraphics{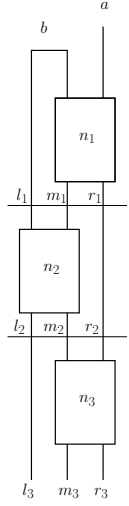}}}
    \caption{$R(n\sb{1}, n\sb{2}, n\sb{3})$ upon numerator closure.}\label{Fi:rn1n2n3}
\end{figure}

We have (see Figure \ref{Fi:rn1n2n3})
\begin{align}\notag
l\sb{3}&=l\sb{2}, \\ \notag
m\sb{3}&=m\sb{2}+(r\sb{2}-m\sb{2})n\sb{3}, \\ \notag
r\sb{3}&=r\sb{2}+(r\sb{2}-m\sb{2})n\sb{3}.
\end{align}

Making use of the results obtained in the preceding subsection,
\begin{align}\notag
l\sb{3}&=a+(b-a)( 1+ n\sb{1}n\sb{2}) ,  \\ \notag
m\sb{3}&=a+(b-a)( 1 + n\sb{1}+n\sb{3}+n\sb{1}n\sb{2}+n\sb{1}n\sb{2}n\sb{3}),
\\ \notag r\sb{3}&=a+(b-a)(n\sb{1}+n\sb{3}+n\sb{1}n\sb{2}n\sb{3}).
\end{align}

We perform a special closure on the tangle in Figure
\ref{Fi:rn1n2n3} and so the coloring system of equations is here
\[
(b-a)(n\sb{1}+n\sb{3}+n\sb{1}n\sb{2}n\sb{3})= 0,
\]
and
\[
\det R(n\sb{1}, n\sb{2}, n\sb{3}) =
n\sb{1}+n\sb{3}+n\sb{1}n\sb{2}n\sb{3}.
\]

\subsubsection{$R(n\sb{1}, n\sb{2}, n\sb{3}, n\sb{4})$}\label{sss:rn1n2n3n4}

\noindent

Analogously,

\begin{align}\notag
l\sb{4}&= a+(b-a)(1+n\sb{1}n\sb{2}+n\sb{1}n\sb{4}+n\sb{3}n\sb{4}
+n\sb{1}n\sb{2}n\sb{3}n\sb{4}), \\ \notag
m\sb{4}&=a+(b-a)(1+n\sb{1}+n\sb{3}+n\sb{1}n\sb{2}+n\sb{1}n\sb{4}+n\sb{3}n\sb{4}
+n\sb{1}n\sb{2}n\sb{3}+n\sb{1}n\sb{2}n\sb{3}n\sb{4}), \\ \notag
r\sb{4}&=a+(b-a)(n\sb{1}+n\sb{3}+n\sb{1}n\sb{2}n\sb{3}).
\end{align}

The coloring system of equations is
\[
(b-a)(1+n\sb{1}n\sb{2}+n\sb{1}n\sb{4}+n\sb{3}n\sb{4}
+n\sb{1}n\sb{2}n\sb{3}n\sb{4}) = 0,
\]
and
\[
\det R(n\sb{1}, n\sb{2}, n\sb{3}, n\sb{4}) =
1+n\sb{1}n\sb{2}+n\sb{1}n\sb{4}+n\sb{3}n\sb{4}
+n\sb{1}n\sb{2}n\sb{3}n\sb{4}.
\]

\subsubsection{$R(n\sb{1}, n\sb{2}, n\sb{3}, n\sb{4}, n\sb{5})$}\label{sss:rn1n2n3n4n5}

\noindent

And analogously,

\begin{align}\notag
l\sb{5}&= a+(b-a)(1+n\sb{1}n\sb{2}+n\sb{1}n\sb{4}+n\sb{3}n\sb{4}
+n\sb{1}n\sb{2}n\sb{3}n\sb{4}), \\ \notag
m\sb{5}&=a+(b-a)(1+n\sb{1}+n\sb{3}+n\sb{5}+n\sb{1}n\sb{2}+n\sb{1}n\sb{4}+n\sb{3}n\sb{4}\\
\notag &\quad
+n\sb{1}n\sb{2}n\sb{3}+n\sb{1}n\sb{2}n\sb{5}+n\sb{1}n\sb{4}n\sb{5}+n\sb{3}n\sb{4}n\sb{5}+
n\sb{1}n\sb{2}n\sb{3}n\sb{4}+n\sb{1}
n\sb{2}n\sb{3}n\sb{4}n\sb{5}),\\ \notag
r\sb{5}&=a+(b-a)(n\sb{1}+n\sb{3}+n\sb{5}
+n\sb{1}n\sb{2}n\sb{3}+n\sb{1}n\sb{2}n\sb{5}+n\sb{1}n\sb{4}n\sb{5}+n\sb{3}n\sb{4}n\sb{5}+n\sb{1}
n\sb{2}n\sb{3}n\sb{4}n\sb{5}).
\end{align}

The coloring system of equations is
\[
(b-a)(n\sb{1}+n\sb{3}+n\sb{5}
+n\sb{1}n\sb{2}n\sb{3}+n\sb{1}n\sb{2}n\sb{5}+n\sb{1}n\sb{4}n\sb{5}+n\sb{3}n\sb{4}n\sb{5}+n\sb{1}
n\sb{2}n\sb{3}n\sb{4}n\sb{5})= 0,
\]
and
\[
\det R(n\sb{1}, n\sb{2}, n\sb{3}, n\sb{4}, n\sb{5})  =
n\sb{1}+n\sb{3}+n\sb{5}
+n\sb{1}n\sb{2}n\sb{3}+n\sb{1}n\sb{2}n\sb{5}+n\sb{1}n\sb{4}n\sb{5}+n\sb{3}n\sb{4}n\sb{5}+n\sb{1}
n\sb{2}n\sb{3}n\sb{4}n\sb{5}.
\]

\bigbreak

\section{The $p\sb{N}$, $p\sb{N}\sp{e}$, and $p\sb{N}\sp{o}$
Polynomials}\label{sect : polynomials}

\noindent



The following definitions help us understanding the structure of the formulas for the $l\sb{N}, m\sb{N},$ and $r\sb{N}$ obtained in Subsection \ref{ss:particularcases}.

\begin{def.}
For each positive integer $N$, let $IEO[N]$ stand for the set of increasing sequences of terms from $\{ 1, 2, \dots , N\}$, whose even terms are even and odd terms are odd, along with the empty sequence $($denoted $\emptyset )$.

If $u\in IEO[N]$, say $u=(u\sb{1}, u\sb{2}, \dots , u\sb{k}) \quad (1\leq k \leq N)$, then let
\[
n\sb{u}:= n\sb{u\sb{1}}n\sb{u\sb{2}} \dots n\sb{u\sb{k}}
\]
If $u=\emptyset$, then let
\[
n\sb{u}:= 1
\]
\end{def.}

\bigbreak

For $N=4$ we obtain, besides the empty sequence, $\emptyset$, the sequences:
\[
(1) \qquad (3) \qquad (1, 2) \qquad (1, 4) \qquad (3, 4) \qquad (1, 2, 3) \qquad (1, 2, 3, 4).
\]

\bigbreak

\begin{def.}\label{def : polys}
For each positive integer $N$, we set:
\[
p\sb{N}:=\sum\sb{u\in IEO[N]}n\sb{u}
\]
\[
p\sb{N}\sp{e}:={\sum\sp{\qquad\prime}\sb{u\in IEO[N]}}  n\sb{u}
\]
where the $'$ denotes deletion of the monomials which are products of an odd number of $n\sb{i}$'s; and
\[
p\sb{N}\sp{o}:=p\sb{N} - p\sb{N}\sp{e}
\]
\end{def.}

\bigbreak

For $N=4$ we obtain,
\[
p\sb{4} = 1 + n\sb{1} + n\sb{3} + n\sb{1}n\sb{2} + n\sb{1}n\sb{4} + n\sb{3}n\sb{4} + n\sb{1}n\sb{2}n\sb{3} + n\sb{1}n\sb{2}n\sb{3}n\sb{4},
\]
\[
p\sb{4}\sp{e} = 1  + n\sb{1}n\sb{2} + n\sb{1}n\sb{4} + n\sb{3}n\sb{4} +  n\sb{1}n\sb{2}n\sb{3}n\sb{4},
\]
\[
p\sb{4}\sp{o} = n\sb{1} + n\sb{3} +  n\sb{1}n\sb{2}n\sb{3},
\]

\bigbreak

Moreover, we recognize that these polynomials were obtained in Subsubsection \ref{sss:rn1n2n3n4} in the expressions of the $l\sb{4}, m\sb{4}$, and $r\sb{4}$. In particular, we can now write:

\begin{align}\notag
l\sb{4}&= a+(b-a)p\sb{4}\sp{e}, \\ \notag
m\sb{4}&=a+(b-a)p\sb{4}, \\ \notag
r\sb{4}&=a+(b-a)p\sb{4}\sp{o}.
\end{align}
with the coloring system of equations:
\[
(b-a)p\sb{4}\sp{e} = 0,
\]
and
\[
\det R(n\sb{1}, n\sb{2}, n\sb{3}, n\sb{4}) =
p\sb{4}\sp{e}.
\]

\bigbreak

Analogous relations are now clear between the polynomials of Definition \ref{def : polys} and other calculations performed in Subsection \ref{ss:particularcases}. We prove below that these relations do not depend on the magnitude of $N$, the number of twists in the rational knots. Before we do that we prove recurrence formulas for the these polynomials which will be useful in the sequel.

\bigbreak

\begin{lem}\label{lem : polyoddeven} Given a positive integer $I$, consider the $2I+1$
variables $n\sb{1}, n\sb{2}, \dots , n\sb{2I}$. We keep the
notation above concerning the $p$ polynomials.
\begin{enumerate}
\item \[
p\sb{2I}=p\sb{2I-1}+n\sb{2I}p\sb{2I-1}\sp{o}.
\]
In particular,
\begin{itemize}
\item \[
p\sb{2I}\sp{o}=p\sb{2I-1}\sp{o} \qquad \text{ and } \qquad
p\sb{2I}\sp{e}=p\sb{2I-1}\sp{e}+n\sb{2I}p\sb{2I-1}\sp{o}.
\]
\end{itemize}

\bigbreak

\item
\[
p\sb{2I+1}=p\sb{2I}+n\sb{2I+1}p\sb{2I}\sp{e}.
\]
In particular,
\begin{itemize}
\item
\[
p\sb{2I+1}\sp{o}=p\sb{2I}\sp{o}+n\sb{2I+1}p\sb{2I}\sp{e}
\qquad \text{ and } \qquad p\sb{2I+1}\sp{e}=p\sb{2I}\sp{e}.
\]
\end{itemize}
\end{enumerate}
 \end{lem} Proof:
 \begin{enumerate}
 \item It is enough to realize that $IEO[2I]$ is the union of $IEO[2I-1]$ with the sequences with an odd number of terms from $IEO[2I-1]$ augmented with the term $2I$.
 \item Omitted since it is analogous to the proof of 1. above.
 \end{enumerate}
 $\hfill \blacksquare $

\bigbreak

\section{The Coloring System of Equations of a Rational Knot: The General Case}\label{sect:colratlinkgen}

\noindent

In this section we prove Theorem \ref{thm:thethm}. In order to do that we first establish the following result.

\begin{figure}[h!]
    \psfrag{...}{\huge $\vdots$}
    \psfrag{2I+1}{\huge $n\sb{2I-1}$}
    \psfrag{2I+2}{\huge $n\sb{2I}$}
    \psfrag{2I+3}{\huge $n\sb{2I+1}$}
    \psfrag{l2I+1}{\huge $l\sb{2I-1}=a+(b-a)p\sb{2I-1}\sp{e}$}
    \psfrag{m2I+1}{\huge $m\sb{2I-1}=a+(b-a)p\sb{2I-1}$}
    \psfrag{r2I+1}{\huge $r\sb{2I-1}=a+(b-a)p\sb{2I-1}\sp{o}=r\sb{2I}$}
    \psfrag{l2I+2}{\huge $l\sb{2I}=l\sb{2I+1}$}
    \psfrag{m2I+2}{\huge $m\sb{2I}$}
    \psfrag{r2I+2}{\huge $r\sb{2I}$}
    \psfrag{l2I+3}{\huge $l\sb{2I+1}$}
    \psfrag{m2I+3}{\huge $m\sb{2I+1}$}
    \psfrag{r2I+3}{\huge $r\sb{2I+1}$}
    \centerline{\scalebox{.50}{\includegraphics{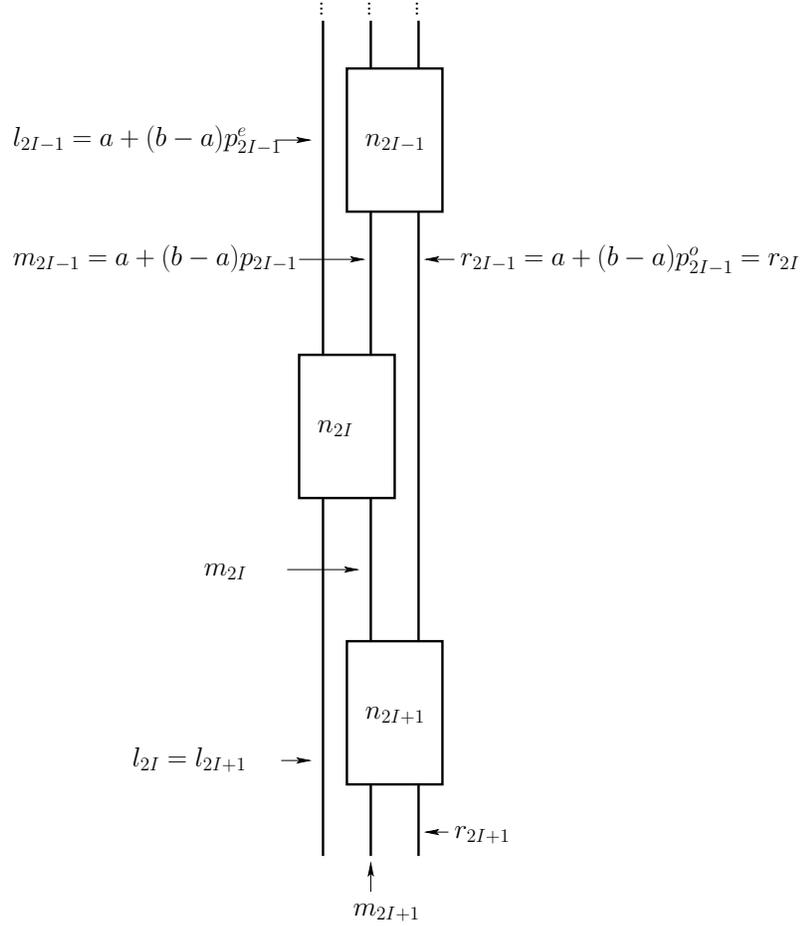}}}
    \caption{(Part of) The induction step in Proposition \ref{prop:gfrat}.}\label{Fi:gfrat}
\end{figure}

\begin{prop}\label{prop:gfrat}
Given a positive integer $N$, fix $N$ integers, $n\sb{1}, n\sb{2},
\dots , n\sb{N}$. Consider the rational knot $R(n\sb{1}, n\sb{2},
\dots , n\sb{N})$ given by the special closure of a $4$-braid. Then,
keeping the notation above, for any $1 \leq i \leq N$,
\begin{align}\notag
l\sb{i}&=a+(b-a)p\sb{i}\sp{e},\\ \notag m\sb{i}&=a+(b-a)p\sb{i}, \\
\notag r\sb{i}&=a+(b-a)p\sb{i}\sp{o}.
\end{align}
\end{prop} Proof: By induction on $N$. The preceding examples show
us that the proposition is true for several particular cases, from
$N=1$ through $N=5$. Assume the statement is true for a given
$N=2I-1$ and consider Figure \ref{Fi:gfrat}. Then,
\begin{align}\notag
l\sb{2I}&=l\sb{2I-1}+n\sb{2I}(m\sb{2I-1}-l\sb{2I-1}) =
a+(b-a)p\sb{2I-1}\sp{e}+n\sb{2I}\biggl( a+(b-a)p\sb{2I-1}-a-
(b-a)p\sb{2I-1}\sp{e}\biggr) \\
\notag & = a+(b-a)p\sb{2I-1}\sp{e}+n\sb{2I} (b-a)\bigl(
p\sb{2I-1}- p\sb{2I-1}\sp{e}\bigr)  =
a+(b-a)p\sb{2I-1}\sp{e}+n\sb{2I} (b-a) p\sb{2I-1}\sp{o} \\
\notag & = a+(b-a)\bigl(
p\sb{2I-1}\sp{e}+n\sb{2I}p\sb{2I-1}\sp{o}\bigr) =
a+(b-a)p\sb{2I}\sp{e}.
\end{align}
The second equality above is a consequence of the induction
hypothesis and the last equality is a consequence of Lemma
\ref{lem : polyoddeven}.

Analogously,
\begin{align}\notag
m\sb{2I}&=m\sb{2I-1}+n\sb{2I}(m\sb{2I-1}-l\sb{2I-1}) =
a+(b-a)p\sb{2I-1}+n\sb{2I}\biggl( a+(b-a)p\sb{2I-1}-a-
(b-a)p\sb{2I-1}\sp{e}\biggr) \\
\notag & = a+(b-a)p\sb{2I-1}+n\sb{2I} (b-a)\bigl( p\sb{2I-1}-
p\sb{2I-1}\sp{e}\bigr)  =
a+(b-a)p\sb{2I-1}+n\sb{2I} (b-a) p\sb{2I-1}\sp{o} \\
\notag & = a+(b-a)\bigl( p\sb{2I-1}+n\sb{2I}p\sb{2I-1}\sp{o}\bigr)
= a+(b-a)p\sb{2I},
\end{align}
and clearly,
\[
r\sb{2I}=r\sb{2I-1}=a+(b-a)p\sb{2I-1}\sp{o}=a+(b-a)p\sb{2I}\sp{o},
\]
using the induction hypothesis and Lemma \ref{lem : polyoddeven}.

In a similar way, assuming now the statement is true for $N=2I$,
we obtain,
\begin{align}\notag
l\sb{2I+1}&= a+(b-a)p\sb{2I+1}\sp{e}, \\ \notag m\sb{2I+1}&=
a+(b-a)p\sb{2I+1}, \\ \notag r\sb{2I+1}&= a+(b-a)p\sb{2I+1}\sp{o}.
\end{align}

This completes the proof. $\hfill \blacksquare $

\bigbreak

Proof (of Theorem \ref{thm:thethm}): We will only prove statement {\it 2.} for {\it 1.} is
similar. As a consequence of Proposition~\ref{prop:gfrat}, the coloring system of equations for
$R(n\sb{1}, n\sb{2}, \dots , n\sb{2I}, n\sb{2I+1})$ reduces to
$a=a+(b-a)p\sb{2I+1}\sp{o}$ which yields at once
\[
(b-a)p\sb{2I+1}\sp{o} = 0 \qquad \text { and } \qquad \det
R(n\sb{1}, n\sb{2}, \dots , n\sb{2I}, n\sb{2I+1}) =
p\sb{2I+1}\sp{o}.
\]
$\hfill \blacksquare $

\section{Counting Spanning Trees of the Graphs Induced by the Checkerboard Shadings of Rational Knots}\label{sect:count}

\noindent

The checkerboard shadings of a knot diagram (Definition
\ref{def:checkandsign}) suggest the definition of a checkerboard
graph in the following way.

\begin{def.}[The checkerboard graph]\label{def:graphcheck}
Given a checkerboard shading of a
knot diagram, we associate to it a graph in the following way. The
vertices of the graph are the shaded regions of the diagram. There
is an edge between two vertices whenever there is a crossing
between the corresponding shaded regions. In particular, two
vertices may be connected by more than one edge.
\end{def.}

\bigbreak

It is a fact that the determinants of alternating knots are given by the counting of the spanning trees of their checkerboard graphs (see \cite{Crowell}). Here we can prove this for rational knots by directly counting the trees and getting the same formula we have derived for rational knots.

\bigbreak

\begin{figure}[h!]
    \psfrag{a1}{\huge $Lou$}
    \psfrag{a2}{\huge $n\sb{2}$}
    \psfrag{...}{\huge $\vdots $}
    \centerline{\scalebox{.50}{\includegraphics{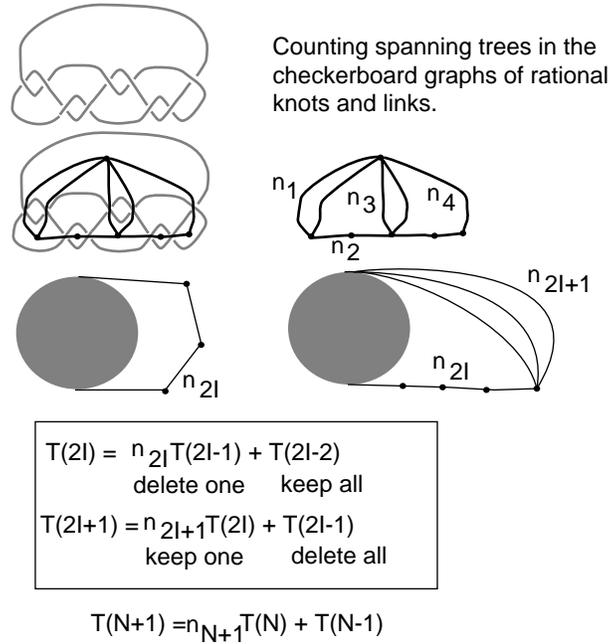}}}
    \caption{Tree count}\label{Fi:Treecount}
\end{figure}

We will assume that an infinite sequence of positive
integers $(n\sb{1}, n\sb{2}, n\sb{3}, \dots )$ has been fixed.
Proposition \ref{prop:recurtree} below yields a recursive relation
for the number of spanning trees as well as a set of instructions
to obtain the new ones from the old ones.

\bigbreak

\begin{prop}\label{prop:recurtree} Letting $T(N)$ stand for the number of spanning trees of
$R(n\sb{1}, n\sb{2}, \dots , n\sb{N})$, then
\[
T(0)=1, \qquad \qquad T(1)=n\sb{1},
\]
\[
T(N+1) = n\sb{N+1}T(N)+T(N-1).
\]
\end{prop}Proof: See Figure \ref{Fi:Treecount}. In this figure we illustrate how the checkerboard graph is obtained from a diagram of
a rational knot, and how the edges in that graph come in either parallel clusters or consecutive sequences depending upon whether the term $n_{N}$ has odd or even index, respectively. With this in mind, one sees that there are two linked formulas for counting the spanning trees in the graph. These formulas are shown in Figure \ref{Fi:Treecount}. In the figure, the instructions, ``delete one'', ``keep all'', ``keep one'', ``delete all'' are appended below the terms in the formulas. These are the instructions for forming the spanning trees from the corresponding sequences and clusters. The reader will find this notation self-explanatory. Once both formulas are written down, we see that they can be summarized by
the single recursion formula of this proposition. This completes the proof.
$\hfill \blacksquare$

\begin{rem} It is a fact that the determinant of an alternating knot or link is equal to the number of spanning trees in its checkerboard graph. Thus we can state as a corollary to this proposition that $T(N) = D_{N}$ where $D_{N}$ denotes the determinant of the rational
knot corresponding to the fixed sequence described above. Thus the formula above gives a recursion formula for the determinants of
rational knots and links. The reader will enjoy verifying that this formula yields our previous formulas for these determinants. As a final remark, note that if all the $n_{N}$ equal $1$, then the sequence $D_{N}$ is the Fibonacci sequence.
\end{rem}

\bigbreak

\bigbreak


\end{document}